\documentclass[12pt]{amsart}
\usepackage[english]{babel}
\usepackage{amsmath}
\usepackage{amsthm}
\usepackage{amssymb}
\usepackage{mathrsfs}
\usepackage{enumerate}
\usepackage[notcite, final, notref]{showkeys}
\usepackage{dsfont}

\usepackage[dvips]{color}
\newtheorem{theorem}{Theorem}[section]

\theoremstyle{definition}
\newtheorem{remark}[theorem]{Remark}

\def\E{\bf E}

\usepackage{xcolor}

\title[]{Sparse Solutions For Inverse Problems in Reproducing Kernel Hilbert Spaces}

\author[T. Qian]{Tao Qian*}
\address{Tao QIAN, Macao Institute of Systems Engineering\\
Macau University of Science and Technology\\
Macau}
\email{tqian@must.edu.mo}
\thanks{*Corresponding author.\\
Funded by The Science and Technology Development Fund, Macau SAR (File no. 0123/2018/A3)}

\def\HH{\mathcal{H}{\mbox -}H_K}
\def\CE{{\bf C}^{\E}}
\def\C{\bf C}

\begin{document}
\maketitle
\begin{abstract} A
linear operator in a Hilbert space defined through the form of Riesz representation naturally introduces a reproducing kernel Hilbert space structure over the range space. Such formulation, called $\mathcal{H}$-$H_K$ formulation in this paper, possesses a built-in mechanism to solve some basic type problems in the formulation by using the basis method, that include identification of the range space, the inversion problem, and the Moore-Penrose pseudo- (generalized) inversion problem. After a quick survey of the existing theory, the aim of the article is to establish connection between this formulation with sparse series representation, and in particular with one called pre-orthogonal adaptive Fourier decomposition (POAFD), the latter being one, most recent and well developed, with great efficiency and wide and deep connections with traditional analysis. Within the matching pursuit methodology the optimality of POAFD is theoretically guaranteed. In practice POAFD offers fast converging numerical solutions.
\end{abstract}

\bigskip


\bigskip

\date{today}
\tableofcontents

\section{Introduction to the $\HH$ formulation, the Basic Problems, and Basis solutions}
\setcounter{equation}{0}
In a Hilbert space if the point-evaluation functional of any point
 is given by the inner product of the function with a function parameterized by the
 point, then we say that the Hilbert space is a reproducing kernel
  Hilbert space (RKHS), and the parameterized function is the
  (unique) reproducing kernel of the RKHS.
  We will start with a formulation of a
linear operator in a general Hilbert space, and lead to a RKHS structure in the range space of the operator. This formulation may be found in a number of sources, and for instance, in \cite{SS2016}. The general Hilbert space is denoted $\mathcal{H}$ with inner product
$\langle\cdot ,\cdot\rangle_{\mathcal H},$ and the linear operator is formulated
with the inner product in the form of the Riesz representation Theorem, as follows. Let $\E$ be an abstract set, usually with a topology.
In our context $\E$ is usually an open set of an Euclidean space, or
 an open set of a domain of one or several complex variables,
 where the elements of $\E$ are treated as parameters.
Associated with each $p\in \E$ there
is an element $h_p\in \mathcal{H}.$ A linear operator $L: \mathcal{H}\to {\CE}$ is defined as
 \begin{eqnarray}\label{operatorL}
 Lf(p)\triangleq\langle f,h_p\rangle_{\mathcal H}.\end{eqnarray}
 where $\CE$ denotes the set of all functions from $\E$ to the complex number field
 ${\bf C}.$ Denote $F(p)=Lf(p).$
Let $N(L)$ be the null space of the operator $L:$
\[ N(L)=\{f\in \mathcal{H}\ |\ L(f)=0\}.\]
$N(L)$ is a closed set in $\mathcal{H}.$ In fact, if
$f_n, f \in \mathcal{H}, f_n\to f$ and $L(f_n)=0,$ then we have
\[ |Lf(p)|=|\langle f-f_n,h_p\rangle_{\mathcal H}|\leq \|f-f_n\|_{\mathcal H}
\|h_p\|_{\mathcal H}\to 0.\]
Thus $Lf=0$ and $f\in N(L).$ As a consequence
we have an orthogonal decomposition for the domain space
\[ \mathcal{H}=N(L)\oplus N(L)^\perp.\]
Accordingly, each $f\in \mathcal{H}$ can be uniquely written as
\[ f=f^-+f^+,\]
where $f^-\in N(L), f^+\in N(L)^\perp.$ We also use the orthogonal projection notations and denote $P_{N(L)^\perp}f=f^+$ and $P_{N(L)}f=f^-,$ where $P_{N(L)^\perp}$ and $P_{N(L)}$ denotes, respectively, the projections
to the closed subspaces $N(L)^\perp$ and $N(L).$ Whenever $F=Lf,$
we have $Lf=Lf^+,$ and $\|f^+\|\leq \|f\|.$ Any solution $g$ for $Lg=F$ has
the form $g=f^++h,$ where $h\in N(L),$ and hence
$\|f^+\|\leq \|g\|.$ Let $R(L)$ denote the range of the operator
$L,$ that is
\[ R(L)=\{F\ |\ \exists f\in \mathcal{H}\ {\rm such \ that\ } F=Lf\}.\] The above particulars
in relation to the orthogonal decomposition of the domain space show that for $F\in R(L)$
the solution $f$ for the equation $Lf=F$ is unique under the minimum norm requirement, and the solution is identical
with $P_{N(L)^\perp}f=f^+.$

We show that $R(L)$ may be equipped with an inner product under which it becomes
a RKHS, denoted $H_K,$ where $K$ stands for
the reproducing kernel. To do this the induced norm of $F=L(f)\in
R(L)$ in the range space
is defined
$$\|F\|_{H_K}\triangleq\|P_{N(L)^\perp}f\|_{\mathcal{H}}.$$ The polarization of the norm
gives rise to an inner product in $R(L)$  denoted
$\langle \cdot,\cdot\rangle_{H_K}.$ The function set $R(L)$ equipped with
the inner product $\langle \cdot,\cdot\rangle_{H_K}$ is named as space
$H_K,$ called the \emph{canonical range space} in relation to $\mathcal{H}$ and $\{h_p\}_{p\in \E}.$
In such way the new Hilbert space $H_K$ is
isometric with $N(L)^\perp$ through the mapping $L.$ Now we show that $K(q,p),$
being being defined as
\[ K(q,p)=\langle h_q,h_p\rangle_{\mathcal{H}},\]
 is
the reproducing kernel of $H_K.$ Alternatively we denote
$K(q,p)=K_q(p).$ We first show
that $h_p\in N(L)^\perp$ and thus $h_p=P_{N(L)^\perp}(h_p).$
For any fixed $p$ in $\E$ the relation $h_p\in N(L)^\perp$ is evidenced by
the fact that for all $f^-\in N(L)$ the relation
\[ 0=L(f^-)(p)=\langle f^-,h_p\rangle_{\mathcal{H}}\]
holds.
Now for $F\in H_K, P_{N(L)^\perp}f=f^+, Lf^+=F, q\in E,$ with
the relation $K_q(p)=\langle h_q,h_p\rangle_{\mathcal{H}}=L(h_q)(p)$, we have 
\begin{eqnarray*}
\langle F,K_q\rangle_{H_K}&=&\langle Lf,L(h_q)\rangle_{H_K}\\
&=&\langle P_{N(L)^\perp}f,P_{N(L)^\perp}h_q\rangle_{\mathcal{H}}\\
&=&\langle f^+,h_q\rangle_{\mathcal{H}}\\
&=& L(f^+)(q)\\
&=& F(q),
\end{eqnarray*}
\def\bD{\bf D}
reproducing the value of the function $F$ at $q\in E.$ In the sequel we will call the above formulation as $\mathcal{H}$-$H_K$ formulation, and $H_K$ the \emph{canonical range space}.

This formulation is as if customized especially for the the complex Hardy spaces: a space having very fundamental impact to
harmonic analysis, complex analysis, as well as to signal analysis.
But it is not: the formulation is a very general and
suitable for all integral, ordinary and partial
differential operators defined in their respective Hilbert spaces (see \cite{SS2016}) in the form of the Riesz representation Theorem.
Below we explain how the complex Hardy space of the unit disc is precisely an example for the $\HH$ formulation.
In the case $\mathcal{H}=L^2(\partial \bD),$
where $\bD$ denotes the complex unit disc and $\partial\bD$
means the boundary of $\bD,$
i.e., the unit circle. $L^2(\partial \bD)$ is facilitated with the inner product
\[ \langle f,g\rangle_{L^2(\partial \bD)}=\frac{1}{2\pi}
\int_0^{2\pi}f(e^{it})\overline{g}(e^{it})dt\]
under which $L^2(\partial \bD)$ is a Hilbert space but itself is not a RKHS.
In the case $\E=\bD.$ For $p\in \bD,$
\[ h_p(e^{it})=\frac{1}{1-\overline{p}e^{it}}\in L^2(\partial \bD).\]
The function $h_p$ is the Szeg\"o kernel of the context
being the Cauchy kernel in the circle arc length measure.
Naturally, for $f\in L^2(\partial\bD),$
$F(p)=\langle f,h_p\rangle_{L^2(\partial \bD)}$ is the Cauchy integral of the boundary
data $f$ over the unit circle. The range space $H_K$ is identical with the complex Hardy
space $H^2(\bD):$
\[ H^2({\bD})=\{F:{\bD} \to {\C}\ |\ F\ {\rm is\ holomorphic\ and}\ \|F\|_{H^2({\bD})}^2=\sup_{0<r<1}\int_0^{2\pi}|F(re^{it})|^2<0\}.\]
A functions $F(z)$ being in this space is equivalent with the condition that $F(z)$ has the Taylor series expansion $F(z)=\sum_{k=0}^\infty c_kz^k$ with $\sum_{k=0}^\infty |c_k|^2<\infty.$ In both the set theoretic and the Hilbert space inner product and norm sense $H^2({\bD})=H_K,$ where functions $F$ in $H_K$ is equipped with the norm $\sum_{k=0}^\infty |c_k|^2.$ We note that
the reproducing kernel of $H^2(\bD)$ is, according to the Cauchy formula,
\[ K(q,p)=K_q(p)=\langle h_q,h_p \rangle_{\mathcal{H}}=\frac{1}{1-\overline{q}p}.\]
The reproducing function of $K_q$ for $F\in H_K$ may be verified through
\[ \langle F,K_q\rangle_{H_K}=\langle f^+,h_q\rangle_{\mathcal{H}}
=\langle f,h_q\rangle_{\mathcal{H}}=F(q).\]
Denote by $H$ the circular Hilbert transform on the circle. The $L^2$ data $f$ on $\partial\bD$ has the decomposition
$f=f^++f^-,$ where $f^+(e^{it})=(1/2)(f+iHf)=\sum_{k=0}^\infty c_ke^{ikt}$ and
$f^-(e^{it})=(1/2)(f-iHf)=\sum_{k=-1}^{-\infty} c_ke^{ikt}.$ $f^\pm$ are also called the analytic signals associated with
$f,$ from the inside and the outside of the disc, respectively. There in particular holds $F(p)=Lf(p)=Lf^+(p)=\langle f^+,h_p\rangle_{L^2(\partial \bD)}.$ The operator $L$ is an isometry mapping between
$f^+$ and $F.$ And, all $f^-,$ non-trivially,
constitute the null space $N(L).$ As a consequence of the Plemelj Theorem
the non-tangential boundary limit of $F$ is identical with $f^+.$ We
note that the inner product of
$H^2(\bD)$ is computed through
the inner product of the isometric subspace $N(L)^\perp$ represented by an integral over the boundary. This is consistent with the $\mathcal{H}$-$H_K$ formulation.
We on the other hand also note that in this Hardy space case there exists an integral with respect to a certain measure
over the whole disc region that
gives rise to the norm as well.
It is referred to as the Littlewood-Paley Identity:  for $F\in H^2(\bD),$
\begin{eqnarray}\label{Hardy} \|F\|^2_{H^2(\bD)}=
|F(0)|^2+2\int_{\bD}|F'(z)|^2\log\frac{1}{|z|}dA(z),\end{eqnarray}
where $dA(z)$ is the normalized area measure of the disc. The polarization of (\ref{Hardy}) gives rise to the integral inner product
 formula of $H_K$ in $\bD$ corresponding to the Littlewood-Paley formula:
\[ \langle F,G\rangle_{H_K}=F(0)\overline{G}(0)+2\int_{\bD}
f'(z)\overline{g}'(z)\log\frac{1}{|z|}dA(z).\]

\bigskip

The $\HH$ formulation is general enough to include a wide class of linear operators including
integral, ordinary and partial differential operators. While the integral operators are obviously included,
we take the differential operators case
 as an illustrative example. In the case the underlying space $\mathcal{H}$ itself is usually a RKHS.
Let $f$ be defined in a RKHS $\tilde{\mathcal{H}}=\tilde{\mathcal{H}}_{\tilde{K}}$ with the reproducing kernel
$\tilde{K}.$ Then
\[ f(x)=\langle f,\tilde{K}_x\rangle_{\tilde{\mathcal{H}}_{\tilde{K}}}.\]
Let $P$ be a  multi-variable-polynomial.
 Then with $\partial=(\partial_1,\cdots,\partial_n)$ we have
 \[ P(\partial)f(x)=\langle f,P(\partial)K_x\rangle_{H_K},\]
 turning the differential operator to an integral operators in a suitable space. In the differential
 operator cases the underlying spaces $\mathcal{H}$ are often Sobolev spaces, being RKHSs, or their subspaces.
    In the $\HH$ formulation
there are three types of questions naturally arising, namely,

\begin{enumerate}[(i)]
  \item How to explicitly represent and
  numerically compute the image function
$F(p)=\langle f,h_p\rangle_{\mathcal{H}}?$

  \item Given a function  $F\in R(L),$ how to represent and numerically approximate the inverse
image function $f$ that satisfies
$F=Lf$ and $\|f\|=\min \{\|g\|\ |\ Lg=F\}?$

  \item Solve the Moore-Penrose psuedo-inverse (generalized inverse) problem:
Assume that the $L$-image space $H_K$ be contained in
a Hilbert space $\tilde{\mathcal{H}}$ as a closed subspace. The question is that for any given function $F\in \tilde{\mathcal{H}},$ find $f\in \mathcal{H}$ such
that $f$ is of the smallest norm in $\mathcal{H}$ and $\|Lf-F\|_{\tilde{\mathcal{H}}}$ is minimized.
\end{enumerate}

There have been studies in relation to these questions (see \cite{SS2016} and its enormous references). There have been ample literature
on reproducing kernel methods in solving various problems of the type of linear operators in Hilbert spaces. Below we summarize what we call as basis method. The basis method as a methodology has existed in the literature. We include here a unified and concise formulation. 

In the $\HH$ formulation $H_K$ is a RKHS. Its normalized reproducing kernels $E_q=\frac{K_q}{\|K_q\|}, q\in \E,$ constitute a dictionary, meaning that the set $\{E_q\ |\ q\in \E\}$ is dense in $H_K.$ The last assertion follows from the reproducing property of the kernels. If the parameter set $\E$ is an open set, and the mapping from $\E$ to the set $\{E_q\ |\ q\in \E\}$ is continuous in the topology of $\mathcal{H},$ then some countable subset $\{E_{q_n}\ |\ q_n\in {\E}, n=1,2,\cdots\},$ can be a complete system of $H_K.$ As a consequence,  $H_K$ contains an orthonormal basis $B_1,B_2,\cdots,$ that is the Gram-Schmidt (G-S) orthonormalization of the collection $\{E_{q_n}\ |\ q_n\in {\E}, n=1,2,\cdots\},$
where
\[ B_n=\frac{E_{q_n}-\sum_{l=1}^{n-1}\langle E_{q_n},B_l\rangle B_l}
{\sqrt{1-\sum_{l=1}^{n-1}|\langle E_{q_n},B_l\rangle|^2}}.\]
We note that in the basis formulation the parameters $q_n, n=1,\cdots,n,\cdots,$ are all distinguished to each other.
Accordingly, we have
\begin{eqnarray}\label{validate}\mathcal{A}_n\mathcal{B}_n=\mathcal{E}_n,\ \quad {\rm and\ thus} \quad \mathcal{B}_n=\mathcal{A}_n^{-1}\mathcal{E}_n,\end{eqnarray}
where for each $n$ the matrix $\mathcal{A}_n$ is of rank $n$ and order $n\times n$ with entries $\langle {E}_i,B_j\rangle_{H_K}, 1\le i,j\le n,$ and the matrices $\mathcal{B}_n$ and $\mathcal{E}_n$ both are of order $n\times 1$ (i.e., column matrices) with entries, respectively, $B_l$ and $E_l, l=1,\cdots, n.$ Due to the property $\langle E_i,B_j\rangle_{H_K}=0$ for all cases $i<j,$ the relations in (\ref{validate}) then be formally extended to the corresponding
infinite matrices as
\begin{eqnarray}\label{valida}\mathcal{A}\mathcal{B}=\mathcal{E}\ \quad {\rm and} \quad \mathcal{B}=\mathcal{A}^{-1}\mathcal{E},\end{eqnarray}
with suitable interpretations of the notations.

To solve the problem (i) one just expands the given $F\in H_K$ into the basis $\{B_l\}_{l=1}^\infty,$ and has
\begin{eqnarray}\label{F1} F={F}_{\mathcal{B}}\mathcal{B},\end{eqnarray}
where $F_{\mathcal{B}}$ is the infinite row matrix consisting of $\langle F,B_l\rangle_{H_K},$ and $\mathcal{B}$ is the infinite column matrix consisting of $B_l.$ Next we solve the inversion problem (ii).
We note that, since $L$ is an isometry from $N(L)^\perp$ to $H_K,$ the inverse operator $L^{-1}$ exists from $H_K$ to $N(L)^\perp,$ being also an isometry.  We have \begin{eqnarray}\label{inverse}L^{-1}F={F}_{\mathcal{B}}L^{-1}\mathcal{B}=
{F}_{\mathcal{B}}\mathcal{A}^{-1}L^{-1}\mathcal{E},\end{eqnarray}
where $L^{-1}\mathcal{E}$ is the infinite column matrix consisting of
the terms $L^{-1}E_{q_n}, n=1,2,\cdots$ The validity of the first equal relation of (\ref{inverse}) is justified by the orthonormality of $L^{-1}\mathcal{B}$ through
a Cauchy sequence argument (also see the proof of Theorem 3.1 below). One can explicitly work out, for any $q\in
\E,$
\[ L^{-1}E_{q}=\frac{L^{-1}K_q}{\|K_q\|_{H_K}}=\frac{h_q}{\|h_q\|_{\mathcal{H}}}.\]

Denoting by $\mathcal{T}$ the infinite column matrix with the entries
$h_{q_n}/\|h_q\|_{\mathcal{H}}, n=1,\cdots,$ we
have
\begin{eqnarray}\label{valid}L^{-1}F=
{F}_{\mathcal{B}}\mathcal{A}^{-1}\mathcal{T}.\end{eqnarray}

Next, we solve the Moore-Penrose pseudo-inverse  problem (iii). The basic assumption is that the space $H_K$ is contained in
a Hilbert space $\tilde{\mathcal{H}}$ as a closed subspace. Let $F$ be the given function in $\tilde{\mathcal{H}}$. The strategy is to expand the projection $G=P_{H_K}(F)$ in $H_K,$ and then expand $G$ into a $\mathcal{B}$-series. Noticing that $F-G$ is perpendicular with $K_q,$ we have
\[ \langle F,K_q\rangle_{\tilde{\mathcal{H}}}=\langle G,K_q\rangle_{\tilde{\mathcal{H}}}=\langle G,K_q\rangle_{H_K}=G(q).\]
Then with
\[ G=\sum_{l=1}^\infty \langle G,B_l\rangle_{H_K}B_l,\]
we have
\begin{eqnarray}\label{MP}
L^{-1}G=\sum_{k=1}^\infty \langle \langle F,K_{\{\cdot\}}\rangle_{\tilde{\mathcal{H}}},B_l\rangle_{H_K}L^{-1}B_l.
\end{eqnarray}
In the matrix notation the above is
\[ L^{-1}G=\{\langle F,K_{\{\cdot\}}\rangle_{\tilde{\mathcal{H}}}\}_{\mathcal{B}}
\mathcal{A}^{-1}L^{-1}\mathcal{E}=
\{\langle F,K_{\{\cdot\}}\rangle_{\tilde{\mathcal{H}}}\}_{\mathcal{B}}
\mathcal{A}^{-1}\mathcal{T},\]
where $\{\langle F,K_{\{\cdot\}}\rangle_{\tilde{\mathcal{H}}}\}_{\mathcal{B}}$ is the row matrix consisting of $\langle \langle F,K_{\{\cdot\}}\rangle_{\tilde{\mathcal{H}}},B_l\rangle_{H_K}, l=1,2,...$ By using the notations $S_1, S_2$ and $S_3$ for the solutions of the problems (i), (ii) and (iii), from (\ref{F1}), (\ref{valid}) and (\ref{MP}), we have

\begin{theorem}\label{Th1}
The solutions to the Problems (i),(ii) and (iii) are respectively given by
\begin{eqnarray}\label{S1} S_1={F}_{\mathcal{B}}\mathcal{B},\end{eqnarray}
 \begin{eqnarray}\label{S2}S_2
={F}_{\mathcal{B}}\mathcal{A}^{-1}\mathcal{T},\end{eqnarray}
and
\begin{eqnarray}\label{S3}
S_3=\{\langle F,K_{\{\cdot\}}\rangle_{\tilde{\mathcal{H}}}\}_{\mathcal{B}}
\mathcal{A}^{-1}\mathcal{T}.
\end{eqnarray}
\end{theorem}

\begin{remark}
The above Problem (iii) is under the assumption that $H_K$ is a subspace of $\tilde{\mathcal{H}}$ that, as a matter of fact, makes a solution straightforward. The example for this is the imbedding of the $L^2$-Bergman space in a complex region into the $L^2$-space in the same region. The more general cases, that is not discussed in the resent paper, include $H_K$ being a set-theoretic subset of $\tilde{\mathcal{H}}$ with a non-isometric imbedding operator $I:\ \|I(f)\|_{\tilde{\mathcal{H}}}\leq \|f\|_{H_K}.$ Such case is, in fact, equivalent in our setting with $\tilde{L}: H_K\to \tilde{\mathcal{H}},$ where $\tilde{L}$ is, in general, a bounded linear operator. Examples for this general cases include, for instance, the imbedding of a Sobolev space into another Sobolev space. 
\end{remark}

We note that the obtained solution formulas are dependent of the basis systems $\mathcal{E},$ $\mathcal{B},$ the transfer matrix $\mathcal{A}.$ They involve complicated computations. The POAFD algorithm proposed in \S 2 is more efficient in computation involving only a limited number of matrices of finite orders for accepted errors.

The rest of the paper introduces a non-basis method, called
pre-orthogonal adaptive Fourier decomposition (POAFD).  
The POAFD method, having been
used in signal and image
analysis, and in system identification, would be, according to the author's knowledge, for the first time introduced to numerical solutions of ODEs, PDEs and integral equations. \S 2 is devoted to the POAFD theory itself.  
 In \S 3 we solve the three types problems by POAFD. The most recent studies show that concrete examples to get numerical solutions using POAFD are all very interesting and significant. As a unified method it is useful whenever the canonical range space $H_K$ is well characterized, or a general kernel $K_q$ is identified. On the other hand, the method itself is helpful to characterize the canonical range space. In the present study we only present the principle of the proposed methods.

\section{POAFD: A Non-Basis Method for Sparse Representation}

Let $H_K$ be the RKHS with kernel function $K(p,q)=K_q(p)=
\langle h_q,h_p\rangle_{\mathcal{H}}$ as in the $\HH$ formulation.  The normalized kernels $E_q=K_q/\|K_q\|_{H_K}, q\in \E,$ constitute a dictionary. Below we will describe the pre-orthogonal adaptive Fourier decomposition (POAFD) algorithm that is available in all Hilbert spaces with a dictionary. Methodology-wise, POAFD belongs to the matching pursuit (or greedy algorithm) type of sparse representations (\cite{MaZ,LT}).  It, however, did not belong to any existing matching pursuit method until it was proposed in \cite{Q2D}. It adopts the idea of Adaptive Fourier Decomposition (AFD) implemented to signals in the classical Hardy spaces. The predecessor AFD was initialized for positive frequency representations of analytic signals, whose algorithm involve the generalized backward shift operator and knowledge of classical Takenaka-Malmquist (TM) system generalizing the Fourier system. It well fits into the frame work of the Beurling-Lax Theorem (\cite{QWa}) and, owing to which, has delicate and deep connections with complex analysis theory, and especially M\"obius transform and Blaschke products. POAFD may be said to be AFD in Hilbert spaces, enhancing delicate analysis due to the fact that it reduces to AFD when underlying Hilbert spaces are replaced by the classical Hardy spaces of one and multiple variables. The AFD algorithm automatically involves multiple parameters (multiple zeros of Blaschke products). Which, in POAFD, corresponds to repeating selections of multiple kernels labelled by the same parameters in the Gram-Schmidt orthogonalization process, when necessary for the optimization principle.  In theoretical development, like in AFD in term of the TM system involving Blaschke products, repeating selections of parameters corresponding to multiple kernels of different levels cannot be avoided. The POAFD maximal selection principle evidences that it is indeed the most effective matching pursuit process. Below we introduce POAFD. To simplify the notation we in the present section borrow the notation
$\{K_q\}_{q\in \E}$ for a collection of functions dense in the underlying Hilbert space, and use $H_K$ for a Hilbert space with a dictionary generalized from $\{K_q\}_{q\in \E}.$ We will not invoke reproducing kernel properties in this section.

For the simplicity, let $\E$ be an open set
 in the complex plane. Let $\{q_1,\cdots,q_n,\cdots,\}$ be an infinite
 sequence of parameters in $\E.$ Denote
\[\tilde{K}_{q_n}=\left(\frac{\partial}{\partial q}\right)^{(l(q_n)-1)}K_q(q_n),\]
where $l(q_n)$ is the number of repeating of the parameter $q_n$ in the $n$-tuple
$\{q_1,\cdots,q_n\}.$ Thus, the function $l(\cdot ),$ as an abuse of notation, is actually dependent of $q_1,\cdots,q_n$ at each test of its value at $q_n.$ We call
 $\tilde{K}_{q_n}, n=1,2,\cdots,$ the
\emph{multiple kernels} that correspond to the parameter sequence in use.  The concept multiple kernel is a necessity
of the pre-orthogonal maximal
 selection principle: Suppose we already have an $(n-1)$-tuple
 $\{q_1,\cdots,q_{(n-1)}\},$ with repetition or not, corresponding to
 the $(n-1)$-tuple $\{\tilde{K}_{q_1},\cdots,\tilde{K}_{q_{n-1}}\}.$
 By doing the G-S orthonormalization process consecutively we obtain
 an equivalent $(n-1)$-orthonormal basis
 $\{B_1,\cdots,B_{n-1}\}.$
 We wish to find a $q_n$ that gives rise to a value being equal, or very close to the following supreme value in the weak-POAFD case:
 \[ \sup\{|\langle G_n,B_n^q\rangle| \ :\ q\in E, q\ne q_1,\cdots,q_{n-1}\}\]
where $G_n$ is the standard remainder
\[ G_n=F-\sum_{k=1}^{n-1}\langle F,B_k\rangle B_k,\] and the finiteness of the supreme is guaranteed by the Cauchy-Schwartz inequality, and
$B_n^q$ be such that $\{B_1,\cdots,B_{n-1},B_n^q\}$ is the G-S orthonormalization
of $\{\tilde{K}_{q_1},\cdots,\tilde{K}_{q_{n-1}}, K_q\},$
given by
 \begin{eqnarray}\label{GS}
 B_n^q
 =\frac{K_q-\sum_{k=1}^{n-1}\langle K_q,B_k\rangle_{H_K}B_k}
 {\sqrt{\|K_q\|^2-\sum_{k=1}^{n-1}|\langle K_q,B_k\rangle_{H_K}|^2}}.\end{eqnarray}
 By definition of
supreme, for any $\rho\in (0,1),$ a parameter $q_n\in \E$
is ready to be found, different from any other previous $q_k, k=1,\cdots, n-1,$  to have
\begin{eqnarray}\label{rou} |\langle G_n,B_n^{q_n}\rangle|\ge \rho
\sup\{\langle G_n,B_n^q\rangle \ :\ q\in \E, q\ne q_1,\cdots,q_{n-1}\}.\end{eqnarray}
The corresponding algorithm for consecutively finding such a sequence  $\{q_n\}_{n=1}^\infty$  is called \emph{Weak-Pre-orthogonal
Adaptive Fourier Decomposition} (Weak-POAFD). With the Weak-POAFD algorithm one
may choose all $q_1,\cdots,q_n$ being distinguished. In many cases, however,
 it happens that
the space satisfies the so called Boundary-Vanishing Condition (BVC):
 For any but fixed $F\in
H_K,$ if $p_n\in \E$ and $p_n\to\partial \E,$ then
\[ \lim_{k\to \infty} |\langle F,E_{p_k}\rangle|=0.\]
If BVC holds, a compact argument leads that
 there exists a point $q_n\in E$
  such that
 \begin{eqnarray}\label{max} |\langle G_n,B_n^{q_n}\rangle|=
\sup\{|\langle G_n,B_n^q\rangle| \ :\ q\in E, q\ne q_1,\cdots,q_{n-1}\}.\end{eqnarray}
When this is the case, the delicate thing is that the limiting point $q_n$ may coincide
with one or several preceding $q_k, k<n.$ In such case it is
the multiple kernel $\tilde{K}_{q_n},$ but not $K_{q_n},$ that has to be used in (\ref{GS}) in doing the G-S process with the preceding $B_1,\cdots,B_{n-1}$ (\cite{Q2D,CQT,qu2018,qu2019}). We note that repeating selection of parameter can be avoided in practice but cannot when doing the theoretical formulation. The theory involving repeating selections in each concrete context is usually not trivial: the beauty of the explicit construction in Szeg\"o kernel and Blaschke products is only a special case, and only for Hardy spaces.  Merely based on the maximal selection principles (\ref{rou}) or (\ref{max}) one can show
\[ F=\sum_{k=1}^\infty \langle F,B_k\rangle_{H_K} B_k\]
(\cite{Q2D,CQT,qu2018,qu2019}).

An order $O(\sqrt{n})$ convergence rate can be proved in a suitably defined subspace
  (\cite{Q2D}). Precisely, for functions
$F$  in the class
\[ \mathcal{M}_M=\{ F\in H_{K}\ |\ \exists \{c_n\} \ {\rm and}\ \ \{K_{q_n}\} \
{\rm such\ that}\ F=\sum_{n=1}^\infty c_nK_{q_n} \ {\rm and}\ \sum_{n=1}^\infty |c_n|\leq M\},\]
the POAFD partial sums satisfy
\[ \|F-\sum_{k=1}^n\langle F,B_k\rangle_{H_K} B_k\|_{H_K}\leq \frac{M}{\sqrt{n}}.\]
We note that POAFD has the same convergence rate as the Shannon expansion into the sinc functions for bandlimited functions.
In the POAFD case the orthonormal system $\{B_1,\cdots,B_n,\cdots\}$ is not necessarily a basis but a system adaptive to the given function $F$ giving rise to fast convergence, being a natural consequence of its maximum selection principle.  It is be just this non-basis violation that gives the capacity of optimal approximation. The algorithm code of POAFD is available at request within the web-page http://www.fst.umac.mo/en/staff/fsttq.html.\\

AFD and POAFD have been seen to have two directions of development. One is $n$-best kernel expansion. That is to determine $n$-parameters at one time, being obviously of better optimality in sparse kernel approximations. $n$-best approximation is motivated by the traditional, yet still open in its ultimate global algorithm: the problem is called the best approximation to Hardy space functions by rational functions of degree not exceeding $n$ (\cite{Bara1986,Baratchart1991,QWM}).
The gradient descending method for cyclic AFD (\cite{QWM}) may be adopted to give $n$-best algorithms in RKHSs. Applications of $n$-best approximation may be found in system identification (\cite{Mi1}), and is usually called model reduction.
The second direction of development of POAFD is related to the Blaschke product and the interpolation in general Hilbert spaces. For existing work along this direction see \cite{ACQS1,ACQS2}. Applications of POAFD, including scalar-valued and matrix-valued Blaschke AFD approximations, to image processing and system identification may be found in
\cite{LZQ,LZQ2,WQLG,WQZS}.

\section{POAFD as Building Block in Solving Problems ({\rm i}), ({\rm ii}) and ({\rm iii})}

POAFD gives the solution of Problem (i) in a fast converging pace and remarkably increases approximation effectiveness. It further makes itself to be the fundamental building block
of the solutions to Problem (ii) and (iii). In this section we come back to the $\HH$ formulation, and, in particular, $H_K$ has reproducing kernel $K.$

\subsection{POAFD Expansion for $F\in H_K:$ the Solution of Problem  ({\rm i})}

Subsequent to what is studied in the last section we have
\begin{eqnarray} S_1={F}_{\mathcal{B}}\mathcal{A}^{-1}\mathcal{E},\end{eqnarray}
where $F_{\mathcal{B}}$ is the infinite row matrix consisting of $\langle F,B_l\rangle_{H_K},$ and $\mathcal{B}$ is the infinite column matrix consisting of $B_l, l=1,2,\cdots,$ being section by section G-S orthonormalization of $\mathcal{E},$ the latter being the infinite column matrix consisting of
the POAFD-selected entries $\tilde{E}_{q_n}, n=1,2,\cdots,$ where $\tilde{E}_{q_n}=
\tilde{K}_{q_n}/\|\tilde{K}_{q_n}\|,$ and $\mathcal{A}$ is the transfer matrix of order $\infty\times \infty$ with entries $\langle \tilde{E}_i,B_j\rangle_{H_K}$ with the property $\langle \tilde{E}_i,B_j\rangle_{H_K}=0$ for $i<j.$

\subsection{The inversion Problem  (ii)}

 The $\HH$ formulation ensures that $L$ is an isometry between $N(L)^\perp$ and
$H_K.$ There, in particular, exists the inverse operator
$L^{-1}$ that maps $F\in H_K$ to the corresponding
$f^+\in N(L)^\perp$ ie., $L^{-1}F=f^+,$ and, in particular,
 $L^{-1}K_q=h_q, q\in E.$
From this, existence and uniqueness of the solution of the
inverse problem follow. Next we work out the explicit series expansion.
 Adaptively expand $F$ by using POAFD:
\begin{eqnarray}\label{Fexpansion} F=\sum_{k=1}^\infty \langle F,B_k\rangle_{H_K}
B_k.\end{eqnarray}
The isometry operator maps the orthonormal system $\{B_k\}_{k=1}^\infty$ to the
orthonormal system $\{L^{-1}B_k\}_{k=1}^\infty.$ We have
\begin{theorem}\label{4.1} With the POAFD-selected parameters $q_1,\cdots,q_n,\cdots,$ there holds
\[ S_2=L^{-1}F=\sum_{k=1}^\infty \langle F,B_k\rangle_{H_K}L^{-1}B_k,\]
where the convergence is in the $\mathcal{H}$-norm sense. In the matrix notation the above solution is written
\begin{eqnarray}\label{inversePOAFD}S_2
={F}_{\mathcal{B}}\mathcal{A}^{-1}\mathcal{T},\end{eqnarray}
where $F_{\mathcal{B}}$, $\mathcal{B}, \mathcal{A}$ and $\mathcal{T}$ are as defined in (\ref{Fexpansion}).

With the $n$-truncated matrices there holds, for $F\in \mathcal{M}_M,$
\begin{eqnarray}\label{except} \| L^{-1}F-{F}_{\mathcal{B}_n}\mathcal{A}_n^{-1}\mathcal{T}_n\|_{\mathcal{H}}\leq \frac{M}{\sqrt{n}}.
 \end{eqnarray}
\end{theorem}
The proof is routine except (\ref{except}). For the self-containing purpose we include the proof for the main convergence part and refer the proof of (\ref{except}) to \cite{Q2D}.
\begin{proof} The $\HH$ formulation shows that there uniquely exists a solution
$f^+=L^{-1}F.$ Since $L^{-1}$ is an isometry between $H_K$ and $N(L)^\perp,$ the system $\{L^{-1}B_k\}$ is orthonormal in the closed subspace $N(L)^\perp.$
Since $\sum_{k=1}^\infty |\langle F,B_k\rangle_{H_K}|^2<\infty,$
the Riesz-Fisher Theorem
concludes that there exists a function $g$ in $N(L)^\perp$ such that
\[ g=\sum_{k=1}^\infty \langle F,B_k\rangle_{H_K}L^{-1}B_k.\]
We need to show that $f^{+}=g.$
It suffices to show
\begin{eqnarray}\label{desired}
 \lim_{n\to \infty}\|f^{+}-\sum_{k=1}^n\langle F,B_k\rangle_{H_K}L^{-1}B_k\|_{\mathcal{H}}^2=0.
 \end{eqnarray}
By using the isometric property of $L^{-1}$ and the relation
(\ref{Fexpansion}), we have
\begin{eqnarray*}
 \lim_{n\to \infty}\|L^{-1}F-L^{-1}(\sum_{k=1}^n\langle F,B_k\rangle_{H_K}B_k)\|^2&=&
 \lim_{n\to \infty}\|F-\sum_{k=1}^n\langle F,B_k\rangle_{H_K}B_k\|^2\\
&=&\lim_{n\to \infty}\|\sum_{k=n+1}^\infty\langle F,B_k\rangle_{H_K}B_k\|^2\\
&=&\lim_{n\to \infty}\sum_{k=n+1}^\infty|\langle F,B_k\rangle_{H_K}|^2=0.
\end{eqnarray*}
The proof is complete.
\end{proof}
This result shows that since we know $L^{-1}K_q=h_q,$ and hence
 $L^{-1}B_k,$ with a POAFD expansion of $F$ we can get
a series expansion solution of the same speed of
convergence for the inverse problem $f^+=L^{-1}F.$

To practically solve an inverse problem under the $\HH$ formulation
the difficulty would be on
finding and characterizing the related objects including $N(L), N(L)^\perp,
K_q.$ In any case $\{h_q\}_{q\in E}$ is a dense subset of $N(L)^\perp.$ In a separate paper we will treat
the special case where the span of
$\{h_q\}_{q\in E}$ is a dense set of $\mathcal{H}$ that corresponds to approximation to identity.

\subsection {The Moore-Penrose Pseudo-Inversion Problem ({\rm iii})}
Problem (iii) is under the assumption that $H_K$
is a closed subspace of a larger Hilbert space
$\tilde{\mathcal{H}}.$  For a given element
${F}\in \tilde{\mathcal{H}}$ the purpose
is to find
\[ {f}\in \mathcal{H} \ {\rm such\ that}\
\|{f}\|_{\mathcal{H}}=\min\{\|\tilde{f}\|_{\mathcal{H}}\ :
\ \|L\tilde{f}-F\|_{\tilde{\mathcal{H}}} \ {\rm is\ minimized}\}.\]

The solution of this problem is divided into two steps.

\noindent{\bf The First Step}
Find the unique function $G\in H_K$ that minimizes $\|F-\tilde{G}\|$
over all $\tilde{G}\in H_K.$
 As given in \S 1, the function $G$ is, in fact,
the projection of $F$ into $H_K,$ denoted $G=P_{H_K}F.$ We already deduced  there $G(q)=\langle F,K_q\rangle_{\tilde{\mathcal{H}}}.$

\noindent{\bf The Second Step} We seek a POAFD series expansion of $G=\langle F,K_q\rangle_{\tilde{\mathcal{H}}}$ as
\[ G=\sum_{k=1}^\infty \langle G,B_k\rangle_{H_K}B_k=\sum_{k=1}^\infty \langle \langle F,K_{(\cdot)}\rangle_{\tilde{\mathcal{H}}},B_k\rangle_{H_K}B_k,\]
where the POAFD is with respect to the reproducing kernel of $H_K$, and the convergence is
 in the $H_K$ norm.
The principle of POAFD shows that the convergence is $O(\frac{1}{\sqrt{n}}).$
We thus have proved
\begin{theorem}\label{MPpseudo}
Under the $\HH$ formulation and the assumption that $H_K$ is a closed subspace
 of $\tilde{\mathcal{H}}$ the solution of the Moore-Penrose pseudo-inverse is
\[S_3=\sum_{k=1}^\infty \langle \langle F,K_{(\cdot)}\rangle_{\tilde{\mathcal{H}}},B_k\rangle_{H_K}B_k,\]
where  $F\in \tilde{\mathcal{H}}$ and the convergence is in the $\tilde{\mathcal{H}}$-norm sense, $K$ is the reproducing kernel of $H_K,$ and
the series is based on the POAFD expansion of the projection function
\[ P_{H_K}F=\langle F,K_{(\cdot)}\rangle_{\tilde{\mathcal{H}}}.\]
By denoting $d_F$ the distance from $F$ to $H_K,$ there holds
\[\| F-\sum_{k=1}^n \langle \langle F,K_{(\cdot)}\rangle_{\tilde{\mathcal{H}}},B_k\rangle_{H_K}B_k\|^2=
d_F^2+O(\frac{1}{n})\]
if $P_{H_K}F\in \mathcal{M}_M.$\end{theorem}

\noindent{\bf Acknowledgement} The author would like to thank Prof Guo, Mao-Zheng, whose discussions improved the author's understanding to the subset imbedding of $H_K$ into a Hilbert space $\tilde{\mathcal{H}}$ which is related to the assumption setting of Problem (iii).

\end{document}